\documentclass[a4paper,11pt]{article}

\usepackage[centertags]{amsmath}
\usepackage[nosort]{cite}
\usepackage[psamsfonts]{amsfonts} 
\usepackage{amsthm,enumerate,amssymb}

\newcommand{\dd}{\mathrm{d}}

\newcommand{\bigtimes}{\mathop{\mbox{\Large $\times$}}} 
\newcommand{\QED}{\hfill{$\square$} \\ \medskip}
\newcommand{\tr}{\operatorname{tr}}
\newcommand{\bbbone}{{\mathchoice {\rm 1\mskip-4mu l} {\rm 1\mskip-4mu l}
{\rm 1\mskip-4.5mu l} {\rm 1\mskip-5mu l}}}
\newcommand{\is}{\mathcal{H}}
\newcommand{\LZ}{\mathcal{L}^\infty_{\mathbb{Z}}}
\newcommand{\Lm}{\mathcal{L}^\infty_{[-t,t]}}
\newcommand{\DZ}{\mathcal{D}_{\mathbb{Z}}}

\theoremstyle{plain}
\newtheorem{theorem}{Theorem}

\begin{document}

\title{Classical dilations \`a la Quantum Probability of Markov evolutions in discrete time}

\author{M.~Gregoratti\thanks{E-mail:
matteo.gregoratti@polimi.it}\\ \\
{\footnotesize \textsl{Dipartimento di Matematica ``F.Brioschi'', Politecnico di
Milano}} \\
{\footnotesize \textsl{Piazza Leonardo da Vinci 32, I-20133 Milano, Italy}}\\ \\
{\small Quaderno di Dipartimento \textbf{QDD 14}}}

\date{}

\maketitle

\begin{abstract}
\noindent We study the Classical Probability analogue of the dilations of a quantum dynamical semigroup 
in Quantum Probability. 
Given a (not necessarily homogeneous) Markov chain in discrete time in a finite state space $E$,  
we introduce a second system, an environment, and a 
deterministic invertible time-homogeneous global evolution of the system $E$ with this environment
such that the original Markov evolution of $E$ can be realized by a proper choice of the initial random 
state of the environment. 
We also compare this dilations with the dilations of a quantum dynamical semigroup in Quantum Probability:
given a classical Markov semigroup, we show that it can be extended to a quantum dynamical 
semigroup for which we can find a quantum dilation to a group of $*$-automorphisms
admitting an invariant abelian subalgebra where this quantum dilation gives just our classical dilation.

\bigskip

\noindent AMS Subject Classification: 60J10, 81S25.
\end{abstract}

\section{Introduction}

We study the analogue in Classical Probability of the dilations in 
Quantum Probability of a quantum dynamical semigroup (QDS) in discrete time. 
A QDS $T^t$ describes the evolution of a quantum system, possibly open, but
``Markovian'', and homogeneous in time. A dilation of $T^t$ by a quantum stochastic flow $j_t$ allows
to represent it by the conditional expectation of a quantum Markov process, analogously to
the representation of a classical Markov semigroup (CMS) by a classical Markov process. 
Anyway, such a dilation in Quantum
Probability enjoys a reacher structure which allows to dilate the semigroup $T^t$ at the same time
also by a unitary group $U^t$, thus showing that the system evolution $T^t$ and the flow
$j_t$ do not contradict the axioms of Quantum Mechanics, i.e.\ that they can arise from a unitary 
time-homogeneous evolution $U^t$ of an isolated bigger system, consisting of the given system and its 
environment \cite{G,K82,K88,M,P,P89}. 

The first aim of this paper is to introduce similar dilations in Classical Probability, choosing a
self-contained approach in a completely classical framework. These classical dilations are interesting, not
only to better understand the relationship between the two probabilistic theories, but also from a simply
classical point of view, to better understand the relationship between Markov processes and 
deterministic invertible homogeneous dynamics.

The second aim of the paper is to show, by appropriate embeddings, that such dilations in Classical
Probability are really analogous to the quantum dilations which inspire them. Even if in Classical Probability
we shall obtain some extra property which has no analogue in Quantum Probability.

We consider here discrete time Markov evolutions in a finite state space. 
A well known theorem due to Birkhoff shows that any doubly stochastic matrix is a convex combination
of permutation matrices. These describe with matrix terminology one-step evolutions which are both
deterministic and invertible. Analogously, a theorem due to A.~S.~Davis \cite{D61}
shows that any stochastic 
matrix, doubly stochastic or not, is a convex combination of deterministic matrices. These describe 
with matrix terminology one-step evolutions which are deterministic but not necessarily invertible. As
a consequence, A.~S.~Davis shows that each finite state Markov chain can be realized as an automaton
with random inputs, thus establishing a connection between Markov chains and automata theory. That is, using a
different terminology, that each finite state Markov chain can be realized via an innovation process.

The relationship between Markov chains and deterministic dynamics can be further analyzed investigating if 
Markov evolutions can be realized as deterministic invertible homogeneous evolutions 
of the system coupled with a second system.
The existence of such representations is theoretically relevant if Markov chains
are applied to phenomena, like physical phenomena for example, for which an underlaying theory
postulates deterministic invertible homogeneous evolutions in absence of noise and external
disturbances. For these phenomena the second system introduced by the dilation models the surrounding
world, the environment, the source of the noise, which is given now a dynamical explanation. 
Of course, here the characterizing property is not simply that the global evolution is deterministic as in
innovation theory, but that it is also invertible and homogeneous in time.

More precisely, we consider a system with finite state space $E$, undergoing a discrete 
time evolution given by a Markov chain, not necessarily homogeneous. 
Then, we introduce an environment  
with its state space $(\Gamma, \mathcal{G})$, a measurable space, together with an
invertible one-step global evolution $\alpha:E\times\Gamma\to E\times\Gamma$. Thus, if $(i,\gamma)$ is the 
state of the
compound system at time 0, then $\alpha^t(i,\gamma)$ is its state at time $t\in\mathbb{N}$,
where hence $\alpha^t$ gives a deterministic invertible homogeneous global evolution. Nevertheless, 
if the environment state is never observed and if initially it is randomly
distributed with some law $Q$ on $(\Gamma,\mathcal{G})$, then the evolution of the observed system turns 
out to be stochastic and, if $\Gamma$, $\mathcal{G}$,
$\alpha$ and $Q$ are properly built, it is given by the original Markov chain.
In this case, we say that 
$(\Gamma,\mathcal{G}, \alpha, Q)$ is a dilation of the Markov evolution in $E$.

Actually in this paper, given only the state space $E$ (arbitrary but finite), we build a \emph{universal} 
dilation
$(\Gamma,\mathcal{G},\alpha,\{Q\})$, where $\{Q\}$ is an entire family of distributions 
which can produce any Markov chain in $E$: 
every Markov chain, homogeneous or not, can be dilated by taking always the same model
$(\Gamma,\mathcal{G},\alpha)$ for the environment and the global evolution,
and by choosing every time the proper distribution $Q$
for the initial state of the environment.
Moreover, our construction, which relays on A.~S.~Davis decomposition of stochastic matrices, allows
to interpret each Markov chain, not only as the stochastic
dynamics resulting from the coupling with an
environment, but at the same time also as an automaton with random inputs, which are now dynamically
provided by the environment itself.

An unpublished result by K\"ummerer \cite{K82} provides a dilation of a homogeneous Markov chain
with a construction similar to ours, but he looks for faithful environment states and he finds an 
interaction $\alpha$ which depends on the chain under consideration. Therefore his construction does not
exhibit that universality which allows us to dilate also non-homogeneous Markov chains.
Our aim is similar also to the aim of Lewis and Maassen \cite{LM84} when they consider classical
mechanics in continuous time and, taken a linear Hamiltonian system modelling a particle and its
environment, they describe how Gibbs states of the whole system lead to stationary Gaussian stochastic
processes for the observables pertaining to the particle under consideration. However, we do not look
for good global states, but for good states $Q$ of the environment alone which lead to Markov
evolutions of the system $E$, our particle, for every independent choice of its initial state.

The paper also proves that our dilations $(\Gamma,\mathcal{G},\alpha,Q)$ are 
restrictions of quantum dilations: every CMS in $E$ admits an extension to a QDS for 
which we can find a quantum dilation which is itself an extension of the dilation of the CMS.
However, we shall not embed the whole universal dilation
$(\Gamma,\mathcal{G},\alpha,\{Q\})$ in the quantum world, as quantum dilations do not exhibit the same
universality and they strictly depend on the QDS under consideration, so that it is not enough
to change the environment state to get another QDS.

In the sequel, given a complex function $f$ on a domain $E$, we shall denote with the same symbol $f$
also its extension on a domain $E\times\Gamma$, $f(i,\gamma)=f(i)$. Similarly, given a map $\phi:E\to E$, we
shall denote with the same symbol $\phi$ also its extension, by tensorizing with the identity, on a domain
$E\times\Gamma$ to $E\times\Gamma$, $\phi(i,\gamma)=\big(\phi(i),\gamma\big)$.

\section{Preliminaries}
We consider a system with finite state space $E=\{1,\ldots,N\}$ and power $\sigma$-algebra $\mathcal{E}$, 
fixed for the whole paper.

\paragraph{Markov evolution.}
We denote by $P=(P_{ij})_{i,j\in E}$ a stochastic matrix on $E$, so that $P_{ij}\geq0$ and $\sum_j
P_{ij}=1$ for every $i$, and we identify the elements of the complex abelian $*$-algebra 
$\mathcal{L}^\infty(\mathcal{E})$, the system random variables
$f:E\to\mathbb{C}$, with the column vectors in $\mathbb{C}^N$, so that every stochastic matrix $P$ in
$E$ defines an operator in $\mathcal{L}^\infty(\mathcal{E})$,
\begin{equation*}
    \big(Pf\big) (i)=\sum_{j\in E}P_{ij}\,f(j),
\end{equation*}
describing a one-step probabilistic evolution. Taken a sequence of stochastic matrices
$(P(t))_{t\geq0}$, $P(0)=\bbbone$, the evolution of a system random variable $f$
from time 0 to time $t\geq0$ is therefore given by
\begin{equation}\label{cobsev0}
    f \mapsto P(1) \cdots P(t) \,f, \qquad \forall f\in \mathcal{L}^\infty(\mathcal{E}).
\end{equation}
We call it Markov evolution and, in the following, we shall denote a sequence 
$(P(t))_{t\geq0}$ simply by $\{P\}$. If the sequence is constant, $P(t)=P$ for every $t$, then
the evolution is homogeneous and it is described by the CMS $(P^t)_{t\geq0}$ in
$\mathcal{L}^\infty(\mathcal{E})$.

\paragraph{Decompositions of stochastic matrices.}
We denote by $D$ a deterministic matrix in $E$, that is a stochastic matrix with
a 1 in each row. Every $D$ describes with matrix terminology a deterministic evolution $\beta$, where
\begin{equation}\label{detmatrix}
    D=(D_{ij})_{i,j\in E}, \qquad \qquad \beta:E\to E, \qquad \qquad D_{ij}=\delta_{\beta(i),j},
\end{equation}
so that $Df=f\circ\beta$.
The invertible (bijective) maps in $E$ correspond to the special cases of permutation matrices.
The deterministic matrices are just the extreme points of the convex set of stochastic matrices
and every $P$ is a convex combination of deterministic matrices,
\begin{equation}\label{dcc}
    P=\sum_{\ell\in L}p_\ell \, D_\ell, \qquad \qquad 
    p_\ell\geq0, \quad \sum_{\ell\in L} p_\ell=1.
\end{equation}
One can find such a decomposition with $N^N$ terms, with the set $L$ labelling all possible
deterministic matrices and weighing each $D_\ell$ with $p_\ell=P_{1\beta_\ell(1)}\cdots 
P_{N\beta_\ell(N)}$. Let us remark that decomposition \eqref{dcc} is not unique and that, for 
any given
$P$, no more than $N^2-N+1$ terms are needed \cite{D61}. Anyway, since we are not going to fix $P$, we
shall employ \eqref{dcc} in the described form: a convex combination of all the deterministic matrices
which can produce any $P$ simply
by changing the weights $(p_\ell)_{\ell\in L}$.

\paragraph{The Markov chain.}
For every Markov evolution $\{P\}$, there exists a Markov chain\\
$\big(\Omega,\mathcal{F},(\mathcal{F}_t)_{t\geq0},(X_t)_{t\geq0},(\mathbb{P}_k)_{k\in E}\big)$ 
with transition probabilities given by $P(t)$, i.e. a discrete time
stochastic process of random variables
$X_t:\Omega\to E$, adapted to a filtration $(\mathcal{F}_t)_{t\geq0}$, and a family of probability
measures $\mathbb{P}_k$, $k\in E$, such that the starting distribution of the process depends on $k$,
$X_0$ has Dirac distribution $\delta_k$ under $\mathbb{P}_k$, but the process always enjoys the 
Markov property with transition 
matrices $P(t)$:
\begin{equation*}
    \mathbb{P}_k(X_{t+1}=j|\mathcal{F}_t) = \mathbb{P}_k(X_{t+1}=j|X_t)
    = P(t+1)_{X_tj}, \qquad \forall k,j\in E,\;t\geq0.
\end{equation*}
Thus, a system random variable $f\in \mathcal{L}^\infty(\mathcal{E})$ has now a stochastic evolution
given by the $*$-unital homomorphism
\begin{equation}\label{stsystev}
    j_t:\mathcal{L}^\infty(\mathcal{E})\to\mathcal{L}^\infty(\mathcal{F}_t), \qquad f\mapsto
    j_t[f]:=f(X_t), \qquad t\geq0,
\end{equation}
and the Markov evolution \eqref{cobsev0} admits the representation
\begin{equation}\label{cobsev1}
    \Big(P(1) \cdots P(t) \,f\Big) (k) = \mathbb{E}_k\big[f(X_t)\big], 
    \qquad \forall f\in \mathcal{L}^\infty(\mathcal{E}).
\end{equation}

\paragraph{Markov chains and automata.}
Let us briefly show a possible realization of the Markov chain 
by means of the decomposition \eqref{dcc}. If $P(t)=\sum_{\ell\in L}p_\ell(t) \, D_\ell$, then one can
take
\begin{gather}
    \nonumber \Omega=E\times L^{\mathbb{N}},\qquad 
    \omega=\big(i,(\ell_n)_{n\in\mathbb{N}}\big),\\
    \nonumber X_0(\omega)=i,\qquad X_t(\omega)=\beta_{\ell_t}\circ\cdots\circ\beta_{\ell_1}(i), \qquad
    Y_t(\omega)=\ell_t,\qquad t\in\mathbb{N}, \\
    \label{Dmc} \mathcal{F}=\sigma(X_0)\otimes\sigma(Y_n;\;n\in\mathbb{N}),
    \qquad \mathcal{F}_t=\sigma(X_0,Y_s;\;1\leq s\leq t),\\
    \nonumber \mathbb{P}_k=\delta_{k}\otimes \big(\bigotimes_{t\in\mathbb{N}}p(t)\big).
\end{gather}
In this way the Markov chain associated to $\{P\}$ is represented as an 
automaton with independent
random inputs:
\begin{equation*}
    X_t=\beta_{Y_t}(X_{t-1}),
\end{equation*} 
so that at every step the set of all possible mappings $\beta_\ell:E\to E$ is available and the
system evolution is determined by the value of the input parameter $\ell$ which is selected randomly
according to $p(t)$ and independently of the previous steps.
If the sequence $P(t)$ is
constant, then the Markov chain is homogeneous and the input parameters can be chosen identically
distributed.
To get a dilation, we shall introduce a bigger $\Omega$, namely $E\times G^\mathbb{Z}$, with a $G$
bigger than $L$ in order to define an invertible one-shot coupling $\phi$, and with $\mathbb{Z}$ 
instead of $\mathbb{N}$ in order to get a dynamics $\alpha^t$ with group properties.

\section{Dilations and universal dilations}

\paragraph{Dilation of a Markov evolution.}
We call dilation of the Markov evolution $\{P\}$ in 
$\mathcal{L}^\infty(\mathcal{E})$ a term
\begin{equation*}
    \Big(\Omega,\mathcal{F},(\mathcal{F}_t)_{t\geq0},(Z_t)_{t\geq0},
    (\mathbb{P}_{k})_{k\in E}\Big)
\end{equation*}
such that
\begin{itemize}
\item every $Z_t=(X_t,\Upsilon_t)$ is a random variable in $(\Omega, \mathcal{F})$ with values in 
$(E\times\Gamma,\mathcal{E}\otimes\mathcal{G})$, being
$(\Gamma,\mathcal{G})$ a fixed measurable space,
\item the term
$\big(\Omega, \mathcal{F}, (\mathcal{F}_t)_{t\geq0}, (X_t)_{t\geq0}, (\mathbb{P}_{k})_{k\in
E}\big)$ is a Markov chain with transition matrices $\{P\}$,
\item the random variable $(X_0,\Upsilon_0)$ has distribution $\delta_k\otimes Q$ under $\mathbb{P}_{k}$,
being $Q$ a fixed distribution on $\mathcal{G}$,
\item there exists an invertible bimeasurable map $\alpha:E\times\Gamma\to E\times\Gamma$ such that
$Z_t=\alpha^t(Z_0)$ for every $t\geq0$.
\end{itemize}

Thus, besides the system $E$, a second system is introduced, an environment with state space
$(\Gamma,\mathcal{G})$. Their states $X_t$ and $\Upsilon_t$ are asked to be random variables on a 
same measurable
space $(\Omega,\mathcal{F})$ such that the global state $Z_t=(X_t,\Upsilon_t)$ undergoes a 
deterministic invertible homogeneous evolution $\alpha^t$. Therefore all the $X_t$ and $\Upsilon_t$ 
are determined by $Z_0$, so that $X_t$ and $\Upsilon_t$ are measurable with respect to
$\sigma(Z_0)=\sigma(X_0,\Upsilon_0)\subseteq\mathcal{F}$ and, depending on the probability 
chosen on $\mathcal{F}$, they are deterministic if and only if $Z_0$ is. Nevertheless, a probability 
$\mathbb{P}_{k}$ typically fixes only the value of $X_0$.
The space $(\Omega, \mathcal{F})$ is also endowed with a filtration $\mathcal{F}_t$. 
Note that only the $X_t$ are asked to be adapted to $\mathcal{F}_t$ so that, in particular,
$\Upsilon_0$ does not have to be $\mathcal{F}_0$-measurable. Therefore
the $X_t$ are not trivially $\mathcal{F}_0$-measurable, even if their values 
are completely determined by the values of $X_0$ and 
$\Upsilon_0$, and,
neglecting the environment, each $\big(\Omega, \mathcal{F}, (\mathcal{F}_t)_{t\geq0}, (X_t)_{t\geq0}, 
\mathbb{P}_{k}\big)$ can be a non trivial stochastic process. What we ask
is that $\big(\Omega, \mathcal{F}, (\mathcal{F}_t)_{t\geq0}, (X_t)_{t\geq0}, 
\mathbb{P}_{k}\big)$ actually is a Markov chain starting from $k$ with transition matrices 
$\{P\}$. At the same time however, this Markov chain is compatible 
with a deterministic, invertible and homogeneous model for the evolution of $E$ coupled with an
environment $\Gamma$. In particular, as $X_0=k$, the whole stochasticity of the process is due only 
to the randomness of the unobserved initial state $\Upsilon_0$ of the environment. 

A dilation gives another interpretation of every evolution \eqref{cobsev0}, 
compatible with \eqref{cobsev1}:
\begin{equation*}
    \Big(P(1) \cdots P(t) \,f\Big)(k) = \mathbb{E}_{k}\big[f(X_t)\big] 
    = \mathbb{E}_{k}\big[f(Z_t)\big] = 
    \mathbb{E}_{k}\Big[f\big(\alpha^t(k,\Upsilon_0)\big)\Big], \qquad \forall f\in
    \mathcal{L}^\infty(\mathcal{E}).
\end{equation*}
Indeed, the stochastic evolution \eqref{stsystev} of a system variable $f\in \mathcal{L}^\infty(\mathcal{E})$ 
is now described by the $*$-unital homomorphism
\begin{equation}\label{stsystev2}
    j_t:\mathcal{L}^\infty(\mathcal{E})\to\mathcal{L}^\infty(\mathcal{F}_t), \qquad     
    f\mapsto j_t[f]:=f(X_t)=f(Z_t)=f\circ\alpha^t(Z_0),
\end{equation}
which is injective as $\alpha^t$ is invertible. And now we could also consider global
random variables $F:E\times\Gamma\to\mathbb{C}$ and their evolution
$F\mapsto F(Z_t)=F\circ\alpha^t(Z_0)$.

\paragraph{Universal dilation.}
Let us denote by $\mathcal{P}$ the set of sequences of stochastic matrices $\{P\}$ in $E$.
We call universal dilation of the Markov evolutions in $\mathcal{L}^\infty(\mathcal{E})$ a term
\begin{equation*}
    \Big(\Omega,\mathcal{F},(\mathcal{F}_t)_{t\geq0},(Z_t)_{t\geq0},
    (\mathbb{P}_{k,\{P\}})_{k\in E,\{P\}\in\mathcal{P}}\Big)
\end{equation*}
such that every $\big(\Omega,\mathcal{F},(\mathcal{F}_t)_{t\geq0},(Z_t)_{t\geq0},
(\mathbb{P}_{k,\{P\}})_{k\in E}\big)$ is a dilation of the corresponding Markov evolution $\{P\}$.
We call universal such a dilation because we ask that the same 
$\Omega$, $\mathcal{F}$, $\mathcal{F}_t$ and $Z_t$ allow
to represent all the Markov evolutions in $\mathcal{L}^\infty(\mathcal{E})$, 
with the change of the probabilities
$\mathbb{P}_{k,\{P\}}$ alone. Therefore, both the environment state space $(\Gamma,\mathcal{G})$ and 
the global evolution
$\alpha$ depend only on the state space $E$, not on the particular Markov evolution to be dilated.

\paragraph{Standard dilation and standard universal dilation.}
In order to show that every state space $E$ admits a universal dilation, now we consider a particular
class of dilations and of universal dilations. 
Let us describe all the special requirements we are interested in for the
dilation of a Markov evolution $\{P\}$.

First of all we want the sample space $\Omega$ to be just $E\times\Gamma$, the state space of the 
global system. As we want it to describe all the possible initial global states, we ask
the random variable $Z_0$ to be the identity function and $X_0$ and $Y_0$ to be the coordinate 
variables: if $\omega=(i,\gamma)$, then $Z_0(\omega)=\omega$, $X_0(\omega)=i$ and
$\Upsilon_0(\omega)=\gamma$. Thus, for all $t\geq0$, $Z_t=(X_t,\Upsilon_t)=\alpha^t\circ
Z_0=Z_0\circ\alpha^t$, $X_t=X_0\circ\alpha^t$ and $\Upsilon_t=\Upsilon_0\circ\alpha^t$. 

We are interested in an
environment $\Gamma=G^\mathbb{Z}=\bigtimes_{n\in\mathbb{Z}} G_n$ with $G_n=G$ 
finite set. In this case the environment state $\gamma\in\Gamma$ has infinitely many components $g_n\in G_n$,
$n\in\mathbb{Z}$, and we introduce also the coordinate variables $Y_n(\omega)=g_n$, the random $n$-th
components of the environment. Then, endowed $G$ with its power $\sigma$-algebra,
we want all the functions so far introduced to be
measurable and so we ask $\mathcal{G}$ to be the natural $\sigma$-algebra 
$\mathcal{G}=\sigma(\Upsilon_0)=\sigma(Y_n;\;n\in\mathbb{Z})$ 
on $\Gamma$, and $\mathcal{F}$ 
to be $\mathcal{E}\otimes\mathcal{G}=\sigma(X_0,\Upsilon_0)$ on $\Omega$. Supposing that at time 0
only $X_0$ is observed and that at each following instant $t$ only the information carried by $Y_t$ is
acquired, we want the filtration $\mathcal{F}_0=\sigma(X_0)$, 
$\mathcal{F}_t=\sigma(X_0, Y_s;\;1\leq s\leq t)$, for $t\geq1$. 

In order to get consistence between these definitions and the global evolution $\alpha$, we ask that
each $Y_t$ is involved in the interaction with the system only once, between time
$t-1$ and time $t$. More precisely, first we ask an
invertible map
\begin{equation}\label{ipc1}
\phi:E\times G\to E\times G,
\end{equation}
denoted by
$\phi(i,g)=\big(\phi^E(i,g),\phi^G(i,g)\big)$, 
giving the one-step coupling between the system and a single environment component. Secondly, we 
introduce the left shift
\begin{equation*}
\vartheta:\Gamma\to \Gamma, \qquad
    (g_n)_{n\in\mathbb{Z}}\in\Gamma,\quad g_n\in G_n \quad \mapsto 
    \quad (g_{n+1})_{n\in\mathbb{Z}}\in\Gamma,\quad g_{n+1}\in G_n.
\end{equation*}
Then, denoted by $\phi_1$ the map $\phi$ in $E\times G_1$, extended 
$\phi_1$ and $\vartheta$ in $\Omega$ by tensorizing with identities, we ask 
\begin{equation*}
\alpha=\vartheta\circ\phi_1.
\end{equation*}
Roughly speaking, when $\alpha$ is applied for the first time between time 0 and time 1, 
the map $\phi_1$ couples the system state $X_0$ with $Y_1$, giving the new system state
$X_1=X_0\circ\alpha=X_0\circ\phi_1=\phi^E(X_0,Y_1)$, and then
the shift $\vartheta$ prepares $Y_2$ for the following interaction 
with $X_1$. So the state of the system at a positive time $t$ is 
$X_t=X_0\circ\alpha^t=\phi^E(X_{t-1},Y_t)$,
which is automatically 
adapted to $\mathcal{F}_t$. Thus $\vartheta$ could be interpreted as a
free evolution of the environment. 
If we explicitly introduce also the random variables $Y^{(t)}_n=Y_n\circ\alpha^t$, the
$n$-th environment components at time $t\geq1$, then
\begin{equation*}Y^{(t)}_n = Y_n\circ\alpha^t = \begin{cases} 
    Y_{n+1}^{(t-1)},& \text{if }n\neq0,\\ 
    \phi^G(X_{t-1},Y_1^{(t-1)}),& \text{if }n=0,
\end{cases} = \begin{cases} 
    Y_{n+t},& \text{if }n\leq-t,\;n\geq1,\\ 
    \phi^G(X_{t-1+n},Y_{n+t}),& \text{if }-t+1\leq n\leq0,
\end{cases}
\end{equation*}
and $X_t=\phi^E(X_{t-1},Y^{(t-1)}_1)$.

At long last, we consider the probabilities $\mathbb{P}_{k}$. Of course, they have to be factorized as
$\delta_k\otimes Q$ on
$\mathcal{F}=\mathcal{E}\otimes\mathcal{G}$. We also require $Q$ to be factorized as
$Q_0\otimes\big(\bigotimes_{n\in\mathbb{N}}q(n)\big)$ on the $\sigma$-algebra
$\sigma(Y_n;\;n\leq0)\otimes\sigma(Y_n;\;n\geq1)$, so that $(Y_n)_{n\leq0}\sim Q_0$ 
while $Y_n\sim q(n)$ for every $n\geq1$ and they are all independent. 
This guarantees the Markov property for the process $X_t$ with respect to $\mathcal{F}_t$
with transition probabilities
\begin{equation}\label{Pq}
P(t)_{ij}=\sum_{g\in G}q_g(t)\,\delta_{\phi^E(i,g),j},\qquad \forall i,j\in E,\;t\geq1,
\end{equation}
so that the only
point is to check if the transition probabilities \eqref{Pq} are the desired ones.

A dilation like this will be called standard in the following.
Summarizing, a dilation is standard if 
\begin{itemize}
\item $\Gamma=G^\mathbb{Z}=\bigtimes_{n\in\mathbb{Z}} G_n$, $\;G_n=G$ finite set with power
$\sigma$-algebra, $\;\mathcal{G}=$ cylindric $\sigma$-algebra,
\item $\Omega=E\times\Gamma, \qquad \omega=(i,\gamma)=(i,(g_n)_{n\in\mathbb{Z}})\in\Omega, 
\quad i\in E, \quad \gamma\in\Gamma, \quad g_n\in G$,
\item $X_0(\omega)=i, \quad\; \Upsilon_0=(Y_n)_{n\in\mathbb{Z}}, \quad\; Y_n(\omega)=g_n,
\quad\; \Upsilon_0(\omega)=\gamma, \quad\; Z_0(\omega)=\omega$,
\item $\mathcal{F}=\mathcal{E}\otimes\mathcal{G}
=\sigma(X_0)\otimes\sigma(Y_n;\;n\leq0)\otimes\sigma(Y_n;\;n\geq1), \qquad 
\mathcal{F}_t=\sigma(X_0, Y_s;\;1\leq s\leq t), \quad t\geq0$,
\item $\mathbb{P}_{k}=\delta_k\otimes Q_0\otimes\big(\bigotimes_{n\in\mathbb{N}}q(n)\big)$,
\item $\alpha=\vartheta\circ\phi_1$, with an invertible $\phi_1=\phi:E\times G_1\to E\times G_1$,
and with the left shift $\vartheta$ on $\Gamma$.
\end{itemize}
Then, for every $t\geq1$ we have
\begin{itemize}
\item $Z_t=\alpha^t\circ Z_0=Z_0\circ\alpha^t=(X_t,\Upsilon_t)$, 
\item $\Upsilon_t=\Upsilon_0\circ\alpha^t=(Y_n^{(t)})_{n\in\mathbb{Z}}$, 
\item $X_t=X_0\circ\alpha^t=\phi^E(X_{t-1},Y_t)=\phi^E(X_{t-1},Y^{(t-1)}_1)$.
\end{itemize}
A standard dilation is therefore specified by the term
$\Big(G,\phi,Q_0\otimes\big(\bigotimes_{n\in\mathbb{N}}q(n)\big)\Big)$. 

With a standard dilation the evolution of every global random variable
is described by the group of $*$-automorphisms $J^t$, $t\in\mathbb{Z}$, where
\begin{equation*} J:\mathcal{L}^\infty(\mathcal{F})\to\mathcal{L}^\infty(\mathcal{F}), \qquad 
J(F) = F\circ\alpha. 
\end{equation*}

Note that a standard dilation always allows to see every Markov chain also as an
automaton with independent random inputs $Y_t$, which now are provided by the environment via the
dynamics $\alpha$. Given a standard dilation 
$\Big(G,\phi,Q_0\otimes\big(\bigotimes_{n\in\mathbb{N}}q(n)\big)\Big)$ for a CMS $P^t$, the random inputs
$Y_t$ can always be chosen independent and identically distributed by taking $q=q(1)$ and replacing 
$Q_0\otimes\big(\bigotimes_{n\in\mathbb{N}}q(n)\big)$ with $q^{\otimes\mathbb{Z}}$.

A universal dilation will be called standard if it is given by a family of standard dilations 
$\big(G,\phi,Q_{\{P\}}\big)$, $\{P\}\in\mathcal{P}$, all of them with the same $G$ and $\phi$. 
Thus we fix the model for the environment and for the global evolution with $G$ and $\phi$, 
and then we require the existence 
of a family of initial distributions $Q_{\{P\}}$ for the environment state, 
each one giving rise to a different Markov evolution for the system. 
A standard universal dilation is therefore specified by
the term $\big(G,\phi,(Q_{\{P\}})_{\{P\}\in\mathcal{P}}\big)$. 
A standard universal dilation is not uniquely determined by the state space
$E$, but it always exists.

\begin{theorem}\label{cud}
For every finite state space $E$, there exists a standard universal dilation\\
$\big(G,\phi,(Q_{\{P\}})_{\{P\}\in\mathcal{P}}\big)$ of the Markov 
evolutions in $\mathcal{L}^\infty(\mathcal{E})$.
\end{theorem}

\noindent {\sl Proof.}
We only have to exhibit a proper set $G$, together with the coupling $\phi$ and the probability 
measures $Q_{\{P\}}$ on $\mathcal{G}$.

Given $E=\{1,\ldots,N\}$ and the set $L$ labelling the all possible maps $\beta:E\to E$, we set
\begin{equation*}
    G=E\times L, \qquad\qquad
    i,j,k\in E, \qquad \ell\in L, \qquad g=(j,\ell)\in G.
\end{equation*}

Arbitrarily fixed $j=1$, we focus on points $(1,\ell)$ in $G$. Thus, taken two points
$\big(i,(1,\ell)\big)\neq\big(i',(1,\ell')\big)$ in $E\times G$, we get
$\big(\beta_\ell(i),(i,\ell)\big)\neq\big(\beta_{\ell'}(i'),(i',\ell')\big)$ 
and so we can find an invertible map
\begin{equation}\label{coupling}
    \phi:E\times G\to E\times G, \qquad\qquad 
    \phi\big(i,(j,\ell)\big) = \begin{cases}\big(\beta_\ell(i),(i,\ell)\big), \quad & \text{if
    }j=1,\\ \quad \ldots \;, & \text{if }j\neq1.\end{cases}
\end{equation}
We choose an arbitrary $\phi$ satisfying \eqref{coupling}.

Given $\{P\}$, 
we fix a decomposition \eqref{dcc} for every $P(t)$, thus obtaining the distributions $p(t)$ on 
$L$; for every distribution $p(t)$ on $L$ we define on $G=E\times L$ the distribution
\begin{equation*}
q(t)=\delta_1\otimes p(t).
\end{equation*}
Thus $\displaystyle P(t)_{ij}=\sum_{\ell\in L}p_\ell(t)\,\delta_{\beta_\ell(i),j}
=\sum_{g\in G}q_g(t)\,\delta_{\phi^E(i,g),j}$ for every $t\in\mathbb{N}$.
Chosen an arbitrary distribution $Q_0$ on $\sigma(Y_n;\;n\leq0)$, if we define
$Q_{\{P\}}=Q_0 \otimes \big(\bigotimes_{t\in\mathbb{N}} q(t)\big)$,
then every stochastic process
$\big(\Omega, \mathcal{F}, (\mathcal{F}_t)_{t\geq0}, (X_t)_{t\geq0}, (\mathbb{P}_{k,\{P\}})_{k\in
E}\big)$ is a Markov chain with transition matrices $\{P\}$, independently of $Q_0$ and 
of the definition of $\phi\big(i,(j,\ell)\big)$ for $j\neq1$. 
Therefore
$\big(G,\phi,(Q_{\{P\}})_{\{P\}\in\mathcal{P}}\big)$
is a standard universal dilation
of the Markov evolutions in $\mathcal{L}^\infty(\mathcal{E})$. 
\QED

In this construction, even if each $g\in G$
has two components, $g=(j,\ell)$, the probability is always
concentrated only  on those $g$ of the kind
$g=(1,\ell)$, but we need the first component $j$ to define an
invertible $\phi$. Analogously, we are considering the evolution only for positive
times so that all the components $g_n$, $n\leq0$, are never involved in the interaction with the
system, but they are needed to define an invertible shift $\vartheta$.

Since the decomposition \eqref{dcc} is not unique, just as the choice of
$Q_0$, we could find other probabilities $\mathbb{P}$ on $\mathcal{F}$ inducing the same 
Markov evolutions for the system. On the other hand, considering non-factorized $Q$ or even 
non-factorized $\mathbb{P}$, we would obtain new stochastic
processes $\big(\Omega, \mathcal{F}, (\mathcal{F}_t)_{t\geq0}, (X_t)_{t\geq0}, \mathbb{P}\big)$, 
without any Markov property guaranteed, which now would depend on the definition of $\phi$ for
$j\neq1$.

Let us remark also that, just because of the universality of the construction, a dilation
$\big(G,\phi,Q_{\{P\}}\big)$ is usually non minimal for a particular 
evolution $\{P\}$.

\paragraph{The cocycle approach to standard dilations.}
Given a Markov evolution $\{P\}$ with a standard dilation 
$\big(G,\phi,Q_0 \otimes \big(\bigotimes_{t\in\mathbb{N}} q(t)\big)\big)$, 
where $G$, $\phi$, $q(t)$ can be defined as in the proof of Theorem \ref{cud} or not, 
as long as we consider only system random variables neglecting the 
environment, we can avoid the shift $\vartheta$ 
and define a deterministic, invertible, but inhomogeneous global evolution 
$\varphi_t$ which never involves the environment
components $g_n$ for $n\leq0$, but which generates the
same Markov chain and the same Markov evolution for the system. Consider indeed
the invertible maps in $\Omega$
\begin{equation*}
    \phi_t:=\vartheta^{-(t-1)} \circ \phi_1 \circ \vartheta^{t-1}, \qquad
    \varphi_t:=\phi_t\circ\cdots\circ\phi_1=\vartheta^{-t}\circ\alpha^t, \qquad
    t\geq1,
\end{equation*}
where $\phi_t$ actually acts only in $E\times G_t$, where it is just the coupling \eqref{coupling},
while $\varphi_t$ actually acts only in $E\times G_1\times\cdots\times G_t$. Then, set
$\varphi_0=\operatorname{Id}$, we get $X_t=X_0\circ\alpha^t=X_0\circ\varphi_t$ for every $t\geq0$. 
Because of the group properties of $\vartheta^t$ 
and $\alpha^t$, the evolution $\varphi_t$ satisfies the cocycle property
\begin{equation*}
\varphi_{t+s}=\vartheta^{-t}\circ\varphi_s\circ\vartheta^t\circ\varphi_t, 
\qquad \forall
t,s\geq1,
\end{equation*}
and, of course, its knowledge is equivalent to the knowledge of $\alpha^t$. 

In particular the $*$-unital injective homomorphism
\eqref{stsystev2}  can be written as
\begin{equation}\label{clfl}
    j_t:\mathcal{L}^\infty(\mathcal{E})\to\mathcal{L}^\infty(\mathcal{F}_t), \qquad f\mapsto
    j_t[f]:=f(X_t)=f\circ\varphi_t, \qquad t\geq0.
\end{equation}
If we denote by $\mathbb{E}_g[f\circ\phi]$ the system random variable in $\mathcal{L}^\infty(\mathcal{E})$
defined by $i\mapsto f\circ\phi(i,g)$, then the stochastic evolution \eqref{clfl} satisfies
\begin{equation}\label{clfleq}
    j_0[f]=f(X_0), \qquad j_t[f]=\sum_{g\in G}
    j_{t-1}\Big[\mathbb{E}_g[f\circ\phi]\Big]\,I_{(Y_t=g)}, \qquad
    \forall f\in\mathcal{L}^\infty(\mathcal{E}), \; t\geq1,
\end{equation}
where $j_t\Big[\mathbb{E}_g[f\circ\phi]\Big]$ are $\mathcal{F}_t$-adapted processes, while the random variables
$I_{(Y_t=g)}$, indicators of the events $(Y_t=g)$, are the increments of the
$\mathcal{F}_t$-adapted processes $N^{g}_t=\sum_{s=1}^tI_{(Y_s=g)}$, $t\geq0$. Hence
Eq.~\eqref{clfleq} can be read as a stochastic equation with respect to the noises $N_t^g$.

We can even reduce the sample space $\Omega$
from $E\times G^\mathbb{Z}$ to $E\times G^\mathbb{N}$, restrict here $\mathcal{F}$,
$\mathcal{F}_t$ and $\mathbb{P}_{k}$, and define $\varphi_t$ in $E\times G^\mathbb{N}$ by \eqref{ipc1}
and \eqref{clfleq}. Anyway, thanks to the cocycle properties of $\varphi_t$, 
it is always possible to introduce later $\bigtimes_{n\leq0}G_n$ and 
the shift $\vartheta$, in order to recover the whole environment
state space $(G^\mathbb{Z},\mathcal{G})$, the evolution $\alpha$ and the initial environment
distribution $Q_0 \otimes \big(\bigotimes_{t\in\mathbb{N}} q(t)\big)$,
so that the two constructions are equivalent and can be considered different descriptions of the same 
situation. 

Choosing the cocycle approach, a standard dilation of a Markov evolution $\{P\}$ is a
Markov chain 
$\big(\Omega, \mathcal{F}, (\mathcal{F}_t)_{t\geq0}, (X_t)_{t\geq0}, (\mathbb{P}_{k,\{P\}})_{k\in E}\big)$,
where
\begin{gather}
    \nonumber\Omega=E\times G^{\mathbb{N}},\qquad
    \omega=\big(i,(g_n)_{n\in\mathbb{N}}\big),\\
    \nonumber X_0(\omega)=i,\qquad Y_t(\omega)=g_t,\qquad 
    X_t=\phi^E(X_{t-1},Y_t), \qquad t\in\mathbb{N},\\
    \label{stMc}\mathcal{F}=\sigma(X_0)\otimes\sigma(Y_t;\;t\in\mathbb{N}),\qquad
    \mathcal{F}_t=\sigma(X_0, Y_s, 1\leq s\leq t), \\
    \nonumber\mathbb{P}_{k}=\delta_{k}\otimes\big(\bigotimes_{t\in\mathbb{N}} q(t)\big).
\end{gather}
This chain is specified by the term 
$\big(G,\phi,\bigotimes_{t\in\mathbb{N}} q(t)\big)$.

Let us underline that the Markov chain \eqref{stMc} is similar to \eqref{Dmc}, 
as also this one is an automaton with independent 
random inputs. Nevertheless, the invertibility of $\phi$ endows this chain 
with a reacher structure because it
implicitly introduces also the
deterministic invertible homogeneous evolution $\alpha^t$.

\section{Dilations of classical Markov semigroups and of quantum dynamical semigroups}

We want to compare a standard dilation with the dilation of a QDS in
Quantum Probability.

Given a Hilbert space $\is$, always complex separable in the paper, 
let us denote its vectors by $h$, or $|h\rangle$ using
Dirac's notation, so that $\langle h'|h\rangle$ denotes the scalar product (linear in $h$) and
$|h'\rangle\langle h|$ denotes the operator $h'' \mapsto \langle h|h''\rangle \,h'$. Let
$\mathcal{B}(\is)$ be the complex $*$-algebra of bounded operators in $\is$. Given a 
Hilbert space $L^2(\mu)$, with a probability measure $\mu$ on some measurable space, and given a measurable
complex function $f$ on the same measurable space,
let $m_f$ denote the multiplication operator
\begin{equation*}
m_f:\operatorname{Dom}(m_f)\to L^2(\mu), \qquad \operatorname{Dom}(m_f)=\{h\in L^2(\mu)\;:\; fh\in L^2(\mu)\},
\qquad m_f\,h=fh,
\end{equation*}
which is bounded if and only if $f\in L^\infty(\mu)$. Let
$\mathcal{D}\big(L^2(\mu)\big)$ be the abelian complex $*$-algebra of bounded multiplication operators in $L^2(\mu)$.
Given a finite space $S$ with its power $\sigma$-algebra $\mathcal{S}$, let $\mu_S$ denote the uniform probability on $(S,\mathcal{S})$ and let $\{|i\rangle\}_{i\in S}$ denote the canonical basis of $L^2(\mu_S)$. Moreover, let us denote $\mathcal{L}^\infty(\mathcal{S})={L}^\infty(\mu_S)$ also by $\mathcal{L}^\infty(S)$.
Given two Hilbert spaces $\is$ and $\mathcal{K}$ and a
trace operator $\tau$ in $\mathcal{K}$, let
$\mathbb{E}_\tau:\mathcal{B}(\is)\otimes\mathcal{B}(\mathcal{K})\to\mathcal{B}(\is)$ denote the conditional
expectation with respect to $\tau$.
Given a vector $\kappa\in\mathcal{K}$, let us denote the conditional
expectation with respect to $|\kappa\rangle\langle\kappa|$ simply by $\mathbb{E}_\kappa$.
In the sequel, given a bounded operator $a$ in $\is$, we shall identify it with its extension
$a\otimes\bbbone_\mathcal{K}$ in $\is\otimes\mathcal{K}$.

\paragraph{Quantum extension of a CMS.}
A discrete-time QDS in $\mathcal{B}(\is)$ is a semigroup $(T^t)_{t\geq0}$ with
$T:\mathcal{B}(\is)\to\mathcal{B}(\is)$ stochastic map, that is a linear, bounded, completely positive, 
normal and identity preserving operator.

In order to extend a CMS in $\mathcal{L}^\infty(\mathcal{E})$ by a QDS in some $\mathcal{B}(\is)$, we take
$\is=L^2(\mu_E)$ and we embed $\mathcal{L}^\infty(\mathcal{E})$ in 
$\mathcal{B}(\is)$ by the $*$-isomorphism $f \mapsto m_f$
between $\mathcal{L}^\infty(\mathcal{E})$ and $\mathcal{D}(\is)$.

We say that a QDS $T^t$ in $\mathcal{B}(\is)$
extends a CMS $P^t$ in $\mathcal{L}^\infty(\mathcal{E})$ if 
\begin{equation*}
    Tm_f = m_{Pf}, \qquad \forall f \in \mathcal{L}^\infty(\mathcal{E}).
\end{equation*}
Such extension always exists. For example, given a stochastic matrix $P$, 
using a representation \eqref{dcc} and notations \eqref{detmatrix}, $P$ is extended by the stochastic map
\begin{equation}\label{Pext}
    Ta = \sum_{\begin{subarray}{c} \scriptscriptstyle\ell\in L \\
    \scriptscriptstyle i\in E \end{subarray}}
    p_\ell\,|i\rangle\langle\beta_\ell(i)|\,a\,|\beta_\ell(i)\rangle\langle i|, \qquad
    \forall a\in\mathcal{B}(\is).
\end{equation}
Note that this extension $T^t$ maps $\mathcal{B}(\is)$ to $\mathcal{D}(\is)$ in only one step. 
The extension of a CMS is not unique at all: for example
every permutation $P$ can be extended also by a $*$-automorphism
$Ta=u^*\,a\,u$, provided that, for every $j$ in $E$, $u\,|j\rangle=P^*\,|j\rangle$ up to a phase 
factor (necessary and sufficient condition).

\paragraph{QP-dilation of a QDS.}
Given a QDS $T^t$ in $\mathcal{B}(\is)$, a typical Quantum Probability construction employs a Hilbert space 
$\mathfrak{Z}$, a unitary operator $V$ in 
$\is\otimes\mathfrak{Z}$ and a unit vector $\upsilon$ in $\mathfrak{Z}$, 
to dilate $T^t$ at the same time by a quantum stochastic flow and 
by a group of $*$-automorphisms. 

Taken infinitely many copies $\mathfrak{Z}_n$ of $\mathfrak{Z}$, the quantum stochastic
flow
\begin{equation*}
    j_t:\mathcal{B}(\is)\to\mathcal{B}(\is)\otimes\mathcal{B}\Big(\bigotimes_{n=1}^t\mathfrak{Z}_n\Big),
    \qquad t\geq0,
\end{equation*}
is the solution of the quantum stochastic equation
\begin{equation*}
    j_0=\operatorname{Id}_{\mathcal{B}(\is)}, \qquad
    j_t(a)=\sum_{zz'}j_{t-1}\Big(\mathbb{E}_{|z\rangle\langle z'|}
    \big[V^*\cdot a\otimes\bbbone_\mathfrak{Z}\cdot V\big]\Big) \otimes
    |z'\rangle\langle z|, \quad t\geq1,
\end{equation*}
where $\{|z\rangle\}$ is a given basis in $\mathfrak{Z}$.

Denoted by $\mathcal{K}$ the infinite tensor product $\bigotimes_{n\in\mathbb{Z}}\mathfrak{Z}_n$ with respect 
to the stabilizing sequence of unit vectors $\psi_n\equiv\upsilon$, denoted by $\Theta$ the left shift 
operator in $\mathcal{K}$, denoted by $V_1$ the operator $V$ in $\is\otimes\mathfrak{Z}_1$, extended in 
$\is\otimes\mathcal{K}$ every operator, consider the unitary operator in $\is\otimes\mathcal{K}$
\begin{equation*}
    U= \Theta\,V_1.
\end{equation*}
The group of $*$-automorphisms is $J^t$, $t\in\mathbb{Z}$, where
\begin{equation*}
J:\mathcal{B}(\is)\otimes\mathcal{B}(\mathcal{K}) \to \mathcal{B}(\is)\otimes\mathcal{B}(\mathcal{K}), 
\qquad J(A)={U^*}\,A\,U.
\end{equation*}
Then $J^t(a)=j_t(a)$ for every $a\in\mathcal{B}(\is)$ and
$t\geq0$.

A term $(\mathfrak{Z},V,\upsilon)$ defines a QP-dilation of the QDS $T^t$ in $\mathcal{B}(\is)$ if
\begin{equation}\label{QPdil}
    T^ta=\mathbb{E}_{\upsilon^{\otimes t}}\big[j_t(a)\big]=
    \mathbb{E}_{\upsilon^{\otimes\mathbb{Z}}}\big[J^t(a)\big], 
    \qquad \forall a\in\mathcal{B}(\is), \;t\geq0.
\end{equation}
The equalities hold for $t\geq0$ if they hold for $t=1$, that is if 
$Ta=\mathbb{E}_{\upsilon}\big[{V^*}\,a\,V\big]$ for all $a$ in 
$\mathcal{B}(\is)$. 

Given $\mathfrak{Z}$, the spatial tensor product $\bigotimes_{n\in\mathbb{Z}}\mathcal{B}(\mathfrak{Z}_n)$ is
a $C^*$-algebra which can be naturally embedded in
$\mathcal{B}(\mathcal{K})$ for every stabilizing sequence $\psi_n$ used in the definition of $\mathcal{K}$. 
Given $V$ and $\upsilon$, the
$C^*$-algebra $\mathcal{B}(\is)\otimes\Big(\bigotimes_{n\in\mathbb{Z}}\mathcal{B}(\mathfrak{Z}_n)\Big)
\subseteq \mathcal{B}(\is)\otimes\mathcal{B}(\mathcal{K})$ is invariant for $J$. Indeed, denoted by 
$\widehat\Theta$ the right shift in $\bigotimes_{n\in\mathbb{Z}}\mathcal{B}(\mathfrak{Z}_n)$, the restriction
of $J$ gives
\begin{equation}\label{qge}
J:\mathcal{B}(\is)\otimes\Big(\bigotimes_{n\in\mathbb{Z}}\mathcal{B}(\mathfrak{Z}_n)\Big) \to
\mathcal{B}(\is)\otimes\Big(\bigotimes_{n\in\mathbb{Z}}\mathcal{B}(\mathfrak{Z}_n)\Big), \qquad
J(A)=V_1^*\,\widehat\Theta(A)\,V_1,
\end{equation}
which is a
$*$-automorphism, is independent of $\upsilon$, and generates a group $J^t$ satisfying \eqref{QPdil}
with $\upsilon^{\otimes\mathbb{Z}}$ pure and locally normal state. 
Thus, at the $C^*$-algebras level, a QP-dilation $(\mathfrak{Z},V,\upsilon)$ of $T^t$ consists of a
$C^*$-algebra $\bigotimes_{n\in\mathbb{Z}}\mathcal{B}(\mathfrak{Z}_n)$ depending on $\mathfrak{Z}$, of a
quantum stochastic flow $j_t$ and of a group of $*$-automorphisms $J^t$ depending on $V$, and of a pure
locally normal state $\upsilon^{\otimes\mathbb{Z}}$ depending on $\upsilon$, 
such that representations \eqref{QPdil} hold.

\paragraph{Quantum extension of a standard dilation.}
Given a stochastic matrix $P$ in $E$ and a standard dilation $\big(G,\phi,q^{\otimes\mathbb{Z}}\big)$ of the CMS $P^t$, we show that $P^t$ can be extended to a QDS $T^t$  
for which we can find a QP-dilation which is itself an extension of $\big(G,\phi,q^{\otimes\mathbb{Z}}\big)$. 
This shows that a standard dilation is a
classical analogue of a QP-dilation.

The first step to study this relationship is to embed the standard dilation
$\big(G,\phi,q^{\otimes\mathbb{Z}}\big)$ in the quantum world. In order to get this embedding at the Hilbert space level, we should introduce a proper measure on $(G^\mathbb{Z},\mathcal{G})$. Taken $\mu_E$ on $E$, the proper measure on $(G^\mathbb{Z},\mathcal{G})$ should give a product measure on 
$(E\times G^\mathbb{Z},\mathcal{E}\otimes\mathcal{G})$ invariant for the deterministic invertible
evolutions $\alpha$, $\vartheta$, $\phi_1$. Thus the natural choice would be the probability measure
$\mu_G^{\otimes\mathbb{Z}}$. Nevertheless, it would be singular with respect to the initial
environment distribution $q^{\otimes\mathbb{Z}}$, so that this latter could not be
obtained from a state in $L^2(\mu_G^{\otimes\mathbb{Z}})$. Therefore now we work at $C^*$-algebras level.

We embed in the quantum world the  spatial tensor product of $C^*$-algebras
\begin{equation*}
\LZ=\mathcal{L}^\infty(\mathcal{E})\otimes\Big(\bigotimes_{n\in\mathbb{Z}}\mathcal{L}^\infty(G_n)\Big).
\end{equation*}
It is the sub $C^*$-algebra of $\mathcal{L}^\infty(\mathcal{F})$ which consists of all continuous
functions in the compact set $E\times G^\mathbb{Z}$, and it is invariant for $\circ\phi_1$, $\circ\vartheta$ and $\circ\alpha$.

Taken again $\is=L^2(\mu_E)$, we introduce the sequence of Hilbert spaces $\mathfrak{Z}_n=L^2(\mu_G)$, $n\in\mathbb{Z}$, the algebras of multiplication operators $\mathcal{D}(\mathfrak{Z}_n)$, and the spatial tensor product
\begin{equation*}
\DZ=\mathcal{D}(\is)\otimes\Big(\bigotimes_{n\in\mathbb{Z}}\mathcal{D}(\mathfrak{Z}_n)\Big),
\end{equation*}
which is an abelian sub $C^*$-algebra of  $\mathcal{B}(\is)\otimes\Big(\bigotimes_{n\in\mathbb{Z}}\mathcal{B}(\mathfrak{Z}_n)\Big)$.

We embed $\LZ$ in $\mathcal{B}(\is)\otimes\Big(\bigotimes_{n\in\mathbb{Z}}\mathcal{B}(\mathfrak{Z}_n)\Big)$ by
the $*$-isomorphism $F \mapsto m_F$ between $\LZ$ and $\DZ$ defined by mapping each
$F$ belonging to 
\begin{equation*}
\Lm:=\mathcal{L}^\infty(\mathcal{E})\otimes\Big(\bigotimes_{n=-t}^t\mathcal{L}^\infty(G_n)\Big), \qquad t\in\mathbb{N},
\end{equation*} 
to the multiplication operator $m_F$ belonging to 
$\mathcal{D}(\is)\otimes\Big(\bigotimes_{n=-t}^t\mathcal{D}(\mathfrak{Z}_n)\Big)
= \mathcal{D}\Big(L^2(\mu_E\otimes\mu_G^{\otimes(2t+1)})\Big)$, abelian subalgebra of
$\mathcal{B}(\is)\otimes\Big(\bigotimes_{n=-t}^t \mathcal{B}(\mathfrak{Z}_n)\Big)
= \mathcal{B}\Big(L^2(\mu_E\otimes\mu_G^{\otimes(2t+1)})\Big) $.

The invertible map $\phi$ in $E\times G$ defines the unitary operator in $\is\otimes\mathfrak{Z}=L^2(\mu_E\otimes\mu_G)$
\begin{equation}\label{psiext}
    V=\sum_{\begin{subarray}{c} \scriptscriptstyle i\in E\\
\scriptscriptstyle g\in G\end{subarray}} |\phi(i,g)\rangle\langle i,g|,
\end{equation}
which is a quantum extension of the deterministic invertible evolution $\phi$ as
\begin{equation*}
V^*\,m_F\,V=m_{F\circ\phi}, \qquad \forall F\in \mathcal{L}^\infty(E)\otimes\mathcal{L}^\infty(G).
\end{equation*}
With respect to the canonical basis of $\mathfrak{Z}$, we have 
$V=\sum_{g,g'\in G}V_{gg'}\otimes|g\rangle\langle g'|$ with
\begin{equation*}
V_{gg'}=\sum_{i,j\in E}|i\rangle\langle i,g|\phi(j,g')\rangle\langle j|.
\end{equation*}

According to Section 3, in the following let $\mathbb{E}_k$ denote the expectation of a random variable on 
$(E\times G^\mathbb{Z},\mathcal{E}\otimes\mathcal{G})$ with respect to the probability measure 
$\delta_k\otimes q^{\otimes\mathbb{Z}}$.
\begin{theorem}\label{cqd} Let $P^t$ be a classical Markov semigroup in a finite state space $E$ and let 
$(G,\phi,q^{\otimes\mathbb{Z}})$ be a standard dilation. 
Let $\is=L^2(\mu_E)$, let $\mathfrak{Z}=L^2(\mu_G)$, let $V$ be the unitary operator \eqref{psiext} in 
$\is\otimes\mathfrak{Z}$, and let $J$ be the $*$-automorphism \eqref{qge} in
$\mathcal{B}(\is)\otimes\Big(\bigotimes_{n\in\mathbb{Z}}\mathcal{B}(\mathfrak{Z}_n)\Big)$. Then
\begin{itemize}
\item[(1)] for every global random variable $F$ belonging to $\LZ$, it holds
\begin{equation}\label{cqd1}
J(m_F)=m_{F\circ\alpha};
\end{equation}
\item[(2)] taken $\upsilon=\sum_{g\in G}\sqrt{q_g}\,|g\rangle\in\mathfrak{Z}$,
\begin{itemize}
\item for every $k\in E$, every $F_1,\ldots,F_n \in\LZ$ and every bounded and continuous function $\eta:\mathbb{C}^n\to\mathbb{C}$, it holds
\begin{equation}\label{cqd2}
    \mathbb{E}_{k}\big[\eta(F_{1},\ldots,F_{n})\big] =  
    \tr\big[\eta(m_{F_1},\ldots,m_{F_n})\,|k\rangle\langle k|\otimes
    |\upsilon^{\otimes\mathbb{Z}}\rangle\langle\upsilon^{\otimes\mathbb{Z}}|\big],
\end{equation}
\item $(\mathfrak{Z},V,\upsilon)$ is a QP-dilation of the quantum dynamical semigroup $T^t$,
\begin{equation}\label{Pext0}
    Ta = \sum_{g\in G} \Bigl(\sum_{g'\in G}\sqrt{q_{g'}}\,V_{gg'}^*\Bigr) \,a\, \Bigl(\sum_{g''\in G}\sqrt{q_{g''}}\,V_{gg''}\Bigr), \qquad    a\in\mathcal{B}(\is),
\end{equation}
which extends the classical Markov semigroup $P^t$.
\end{itemize}
\end{itemize}
\end{theorem}

\noindent {\sl Proof.}
(1) For every $F$ belonging to an algebra $\Lm$,
\begin{equation*}
    \widehat\Theta(m_F)=m_{F\circ\vartheta}, \qquad
    V_1^*\,m_F\,V_1=m_{F\circ\phi_1}, \qquad
    J(m_F)=m_{F\circ\alpha},
\end{equation*}
so that, by continuity, Eq.~\eqref{cqd1} is proved for all $F\in\LZ$.

(2) For every $t\in\mathbb{N}$, the pure locally normal state 
$|\upsilon^{\otimes\mathbb{Z}}\rangle\langle\upsilon^{\otimes\mathbb{Z}}|$ restricted to 
$\bigotimes_{n=-t}^t\mathcal{B}(\mathfrak{Z}_n)$ gives the state 
$|\upsilon^{\otimes(2t+1)}\rangle\langle\upsilon^{\otimes(2t+1)}|$, where the unit vector 
$\upsilon^{\otimes(2t+1)}$ belongs to $ \bigotimes_{n=-t}^t\mathfrak{Z}_n=L^2(\mu_G^{\otimes(2t+1)})$ and has
$|\upsilon^{\otimes(2t+1)}|^2=\dd q^{\otimes(2t+1)}/\dd\mu_G^{\otimes(2t+1)}$.
Thus 
\begin{equation*}
    \mathbb{E}_{k}[F] =  
    \tr\Big[m_{F}\cdot|k\rangle\langle k|\otimes|\upsilon^{\otimes\mathbb{Z}}\rangle
    \langle\upsilon^{\otimes\mathbb{Z}}|\Big]
\end{equation*}
for all $F\in\Lm$ and hence, by continuity, for all $F\in\LZ$. Then Eq.~\eqref{cqd2} is proved, 
as $\eta(F_{1},\ldots,F_{n})$ belongs to $\LZ$
whenever $F_1,\ldots,F_n\in\LZ$ and $\eta$ is bounded and continuous.

Moreover, for every $a\in\mathcal{B}(\is)$,
\begin{multline*}
\mathbb{E}_{\upsilon}\big[V^*\,a\,V\big]= \sum_{j',j''\in E}|j'\rangle\,\langle j'|\otimes\langle\upsilon|\,V^*\,a\,V\, |j''\rangle\otimes|\upsilon\rangle\,\langle j''| \\
= \sum_{\begin{subarray}{c} \scriptscriptstyle j',j''\in E \\ \scriptscriptstyle g',g''\in G \end{subarray}} |j'\rangle\,\sqrt{q_{g'}}\,\langle j',g'|\,V^*\,a\,V\, |j'',g''\rangle\,\sqrt{q_{g''}}\,\langle j''|
= \sum_{\begin{subarray}{c} \scriptscriptstyle j',j''\in E \\ \scriptscriptstyle g',g''\in G \end{subarray}} \sqrt{q_{g'}q_{g''}}\,|j'\rangle\langle \phi(j',g')|\,a\,|\phi(j'',g'')\rangle\langle j''| \\
= \sum_{g,g',g''\in G} \sqrt{q_{g'}q_{g''}}\,V_{gg'}^*\,a\,V_{gg''}
= Ta. 
\end{multline*}
Then
\begin{equation*}
    T m_f = \sum_{\begin{subarray}{c} \scriptscriptstyle g\in G \\ \scriptscriptstyle i\in E \end{subarray}}
    q_g\,|i\rangle\langle\phi^E(i,g)|\,m_f\,|\phi^E(i,g)\rangle\langle i| = m_{Pf}, \qquad \forall f\in \mathcal{L}^\infty(\mathcal{E}).
\end{equation*}
\QED

Equation \eqref{cqd1} shows that the group of $*$-automorphisms $J^t$ defined by $(\mathfrak{Z},V)$ is a 
quantum extension of the deterministic invertible homogeneous evolution $\alpha^t$ defined by $(G,\phi)$. 
Then, automatically, also the quantum stochastic flow $j_t$ extends the classical $j_t$. Moreover, Equation \eqref{cqd2} shows that the whole standard dilation $(G,\phi,q^{\otimes\mathbb{Z}})$ of $P^t$ is extended by the QP-dilation $(\mathfrak{Z},V,\upsilon)$ of $T^t$ as the states $q^{\otimes\mathbb{Z}}$ and $\upsilon^{\otimes\mathbb{Z}}$ 
of the environments, together with the dynamics $\circ\alpha$ and
$J$, give rise to the same joint distribution for the trajectories of the global random variables 
$F$ belonging to $\LZ$ and for the trajectories of the corresponding normal operators
$m_F$ belonging to $\mathcal{B}(\is)\otimes\Big(\bigotimes_{n\in\mathbb{Z}}\mathcal{B}(\mathfrak{Z}_n)\Big)$.
This happens for every possible starting state $k$ of the classical Markov system and for the
corresponding state $|k\rangle\langle k|$ of its quantum counterpart.

Let us note also that if $(G,\phi,q^{\otimes\mathbb{Z}})$ is built as in Theorem \ref{cud}, then
\begin{equation*}
V_{(i,\ell)(1,\ell')}=\delta_{\ell\ell'}\,|\beta_\ell(i)\rangle\langle i|,
\end{equation*}
and the stochastic map \eqref{Pext0} becomes just the stochastic map \eqref{Pext}.

\end{document}